# LIMIT DISTRIBUTION THEORY FOR MAXIMUM LIKELIHOOD ESTIMATION OF A LOG-CONCAVE DENSITY

BY FADOUA BALABDAOUI, KASPAR RUFIBACH[1] AND JON A. WELLNER[2]

*Universite Paris-Dauphine and University of Göttingen, University of Zurich and University of Washington*

We find limiting distributions of the nonparametric maximum likelihood estimator (MLE) of a log-concave density, that is, a density of the form $f_0 = \exp \varphi_0$ where $\varphi_0$ is a concave function on $\mathbb{R}$. The pointwise limiting distributions depend on the second and third derivatives at 0 of $H_k$, the "lower invelope" of an integrated Brownian motion process minus a drift term depending on the number of vanishing derivatives of $\varphi_0 = \log f_0$ at the point of interest. We also establish the limiting distribution of the resulting estimator of the mode $M(f_0)$ and establish a new local asymptotic minimax lower bound which shows the optimality of our mode estimator in terms of both rate of convergence and dependence of constants on population values.

## 1. Introduction.

1.1. *Log-concave densities.* A probability density $f$ on the real line is called log-concave if it can be written as

$$f(x) = \exp \varphi(x)$$

for some concave function $\varphi : \mathbb{R} \to [-\infty, \infty)$. We let $\mathcal{LC}$ denote the class of all log-concave densities on $\mathbb{R}$. As shown by Ibragimov (1956), a density function $f$ is log-concave if and only if its convolution with any unimodal density

Received January 2008.
[1]Supported in part by the Swiss National Science Foundation.
[2]Supported in part by NSF Grant DMS-05-03822 and NI-AID Grant 2R01 AI291968-04.

is again unimodal. Thus, the class of log-concave densities is often referred to as the class of "strongly unimodal" densities. Furthermore, the class $\mathcal{LC}$ of log-concave densities is exactly the class of Polyá frequency functions of order 2, $PFF_2$ as noted by Pal, Woodroofe and Meyer (2007); see also Dharmadhikari and Joag-Dev (1988), page 150, and Marshall and Olkin (1979), page 492.

The log-concave shape constraint is appealing for many reasons:

(1) Many parametric models, for a certain range of their parameters, are in fact log-concave, for example, normal, uniform, gamma$(r, \lambda)$ for $r \geq 1$, beta$(a, b)$ for $a \geq 1$ and $b \geq 1$, generalized Pareto, Gumbel, Fréchet, logistic or Laplace, to mention only some of these models. Therefore, assuming log-concavity offers a flexible nonparametric alternative to purely parametric models. Note that a log-concave density need not be symmetric.

(2) Every log-concave density is automatically unimodal. Furthermore, log-concavity of a density $f$ immediately implies specific shape constraints for certain functions derived from $f$ [see Barlow and Proschan (1975), Marshall and Olkin (1979, 2007), Dharmadhikari and Joag-Dev (1988), An (1998) and Bagnoli and Bergstrom (2005)]. Thus, having an estimator (and its limiting distribution) for $f$ at hand provides, almost automatically, estimators (and limiting distributions) for those functions. Corollary 2.3 illustrates this for the hazard rate.

(3) Although the nonparametric MLE of a unimodal density does not exist [see, e.g., Birgé (1997)], the nonparametric MLE of a log-concave density exists, is unique and has desirable consistency and rates of convergence properties. Thus, the class of log-concave (or strongly unimodal) densities may be a useful and valuable surrogate for the larger class $\mathcal{U}$ of unimodal densities.

(4) Tests for multimodality and mixing can be based on a semiparametric model with densities of the form $f_{c,\varphi}(x) = \exp(\varphi(x) + cx^2)$, where $\varphi$ is concave and $c > 0$, as shown by Walther (2002).

(5) Chang and Walther (2007) further show that the EM-algorithm can be extended to work for log-concave component densities.

(6) First attempts to estimate a log-concave density in $\mathbb{R}^d$ were made by Cule, Gramacy and Samworth (2007).

(7) The log-concave density estimator can be used to improve accuracy in the estimation of the so-called "tail index" of a generalized Pareto distribution [see Müller and Rufibach (2009)].

(8) It should be noted that no arbitrary choices such as bandwidth, kernel or prior are involved in the estimation of a log-concave density; these are all obviated by this shape restriction.

(9) We expect good adaptivity properties of the MLE $\widehat{f}_n$ in the class $\mathcal{LC}$.

For properties of (random variables with) log-concave densities, we refer to Dharmadhikari and Joag-Dev (1988), Marshall and Olkin (1979) and Rufibach (2006). Log-concavity of a density $f$ implies certain shape constraints for functions derived from $f$, such as the distribution function, the tail or hazard function. See An (1998) for comparisons with the related notion of a log-convex density.

1.2. *Log-concave density estimation.* Now let $X_{(1)} < X_{(2)} < \cdots < X_{(n)}$ be the order statistics of $n$ independent random variables $X_1, \ldots, X_n$, distributed according to a log-concave probability density $f_0 = \exp \varphi_0$ on $\mathbb{R}$. The distribution function corresponding to $f_0$ is denoted by $F_0$.

The maximum likelihood estimator (MLE) of a log-concave density was introduced in Rufibach (2006) and Dümbgen and Rufibach (2009). Algorithmic aspects were treated in Rufibach (2007) and in a more general framework in Dümbgen, Hüsler and Rufibach (2007), while consistency with respect to the Hellinger metric was established by Pal, Woodroofe and Meyer (2007), and rates of convergence of $\widehat{f}_n$ and $\widehat{F}_n$ were established by Dümbgen and Rufibach (2009). Since the derivation of the MLE of a log-concave density is extensively treated in these references, we only briefly recall its definition and the properties relevant for this paper.

If $\mathcal{C}$ denotes the class of all concave functions $\varphi : \mathbb{R} \to [-\infty, \infty)$, the estimator $\widehat{\varphi}_n$ of $\varphi_0$ is the maximizer of the "adjusted" criterion function

$$L(\varphi) = \int_{\mathbb{R}} \varphi(x) \, d\mathbb{F}_n(x) - \int_{\mathbb{R}} \exp \varphi(x) \, dx$$

over $\mathcal{C}$, where $\mathbb{F}_n$ is the empirical distribution function of the observations. The log-concave density estimator is then $\widehat{f}_n := \exp \widehat{\varphi}_n$, which exists and is unique.

1.3. *Characterization of $\widehat{\varphi}_n$.* For any continuous piecewise linear function $h_n : [X_{(1)}, X_{(n)}] \to \mathbb{R}$, such that the knots of $h_n$ coincide with (some of) the order statistics $X_{(1)}, \ldots, X_{(n)}$, introduce the set of knots $\widehat{\mathcal{S}}_n(h_n)$ of $h_n$ as

$$\widehat{\mathcal{S}}_n(h_n) := \{t \in (X_{(1)}, X_{(n)}) : h_n'(t-) > h_n'(t+)\} \cup \{X_{(1)}, X_{(n)}\}.$$

Dümbgen and Rufibach (2009) found that $\widehat{\varphi}_n$ is piecewise linear, that $\widehat{\varphi}_n = -\infty$ on $\mathbb{R} \setminus [X_{(1)}, X_{(n)}]$ and that the knots of $\widehat{\varphi}_n$ only occur at (some of the) ordered observations $X_{(1)} < \cdots < X_{(n)}$. The latter property is entirely different from the estimation of a $k$-monotone density for $k > 1$ (see below), where the knots fall strictly between observations with probability equal to 1.

According to Theorem 2.4 in Dümbgen and Rufibach (2009), the estimator $\widehat{\varphi}_n$ has the following characterization. For $x \geq X_{(1)}$ (recall that $\widehat{\varphi}_n := -\infty$ outside $[X_{(1)}, X_{(n)}]$), define the processes

$$\widehat{F}_n(x) := \int_{X_{(1)}}^{x} \exp(\widehat{\varphi}_n(t)) \, dt, \qquad \widehat{H}_n(x) := \int_{X_{(1)}}^{x} \widehat{F}_n(t) \, dt,$$

$$\mathbb{H}_n(x) := \int_{X_{(1)}}^{x} \mathbb{F}_n(t)\,dt = \int_{-\infty}^{x} \mathbb{F}_n(t)\,dt.$$

Then, the concave function $\widehat{\varphi}_n$ is the MLE of the log-density $\varphi_0$ if, and only if,

(1.1)	$\widehat{H}_n(x) \begin{cases} \leq \mathbb{H}_n(x), & \text{for all } x \geq X_{(1)}, \\ = \mathbb{H}_n(x), & \text{if } x \in \widehat{\mathcal{S}}_n(\widehat{\varphi}_n). \end{cases}$

1.4. *Other shape constraints.* Maximum likelihood estimation of a monotone density $f_0$ on $[0,\infty)$ was first studied by Grenander (1956). Under the assumption that $f_0$ is $C^1$ in a neighborhood of a point $x_0 > 0$, such that $f_0'(x_0) < 0$, Prakasa Rao (1969) established the (local) asymptotic distribution theory of the Grenander estimator $\hat{f}_n$:

$$n^{1/3}(\hat{f}_n(x_0) - f_0(x_0)) \xrightarrow{d} |f_0'(x_0)f_0(x_0)/2|^{1/3}\mathbb{Z},$$

where $\mathbb{Z}$ is the slope at zero of the (least) concave majorant of the process $W(t) - t^2$, $t \in \mathbb{R}$ for two-sided Brownian motion $W$ starting at 0.

Under the assumption that the true density $f_0$ is convex on $[0,\infty)$ and that $f_0$ is $C^2$ in a neighborhood of $x_0$ with $f_0''(x_0) > 0$, Groeneboom, Jongbloed and Wellner (2001b) show that the MLE $\hat{f}_n$ (as well as the least squares estimator of $f_0$) satisfies

$$n^{2/5}(\hat{f}_n(x_0) - f_0(x_0)) \xrightarrow{d} (24^{-1}f_0^2(x_0)f_0''(x_0))^{1/5}\mathbb{H}''(0),$$

where $\mathbb{H}$ is a particular upper invelope of an integrated two-sided Brownian motion $+t^4$ [see also Groeneboom, Jongbloed and Wellner (2001a)].

The classes of monotone and convex decreasing densities are particular cases of the class of $k$-monotone densities. Modulo a spline interpolation conjecture, Balabdaoui and Wellner (2007) were able to adapt the approach of Groeneboom, Jongbloed and Wellner (2001b) to this general class of densities.

We find that log-concave estimation shares many similarities with the aforementioned shape-constrained estimation problems. In particular, the limiting distribution of the MLE, our nonparametric estimator, involves a stochastic process whose second derivative is concave and which stays below an integrated Brownian motion minus $t^{k+2}$. The even integer $k$ determines the number of vanishing derivatives of the true concave function $\varphi_0$ at the estimation point $x_0$. Using Theorem 2.1, one can derive a procedure for estimation of $k$. This is relevant in practical applications of our results, that is, construction of confidence intervals for the mode using the limiting distribution given in Theorem 2.1. These problems are the subject of ongoing research.

1.5. *Organization of the paper.* In Section 2, we establish the limiting distributions of the ML estimators, $\widehat{\varphi}_n$ and $\widehat{f}_n$, at a fixed point $x_0 \in \mathbb{R}$ under some specified working assumptions. The characterization of either $\widehat{\varphi}_n$ or $\widehat{f}_n$ given in (1.1) coincides, except for the direction of the inequality, with that of the least-squares estimator of a convex decreasing density, studied by Groeneboom, Jongbloed and Wellner (2001b); see their Lemma 2.2, page 1657. This enables us to adopt the general scheme of the proof in their paper.

Log-concave densities $f$ and their logarithm $\varphi$ can easily have vanishing second and higher derivatives at fixed points; an explicit example will be given in Section 2. Thus, the formulation of our asymptotic results allows higher derivatives of the concave function $\varphi_0$ to vanish at the estimation point. This is somewhat more general than the assumptions of Groeneboom, Jongbloed and Wellner (2001b) (where a natural assumption is that the second derivative is positive at the point of interest, but similar vanishing of second derivatives and existence of a nonzero higher order derivative can also easily occur), but it is analogous to the results of Wright (1981) and Leurgans (1982) for nonparametric estimation of a monotone regression function. Similar results for the Grenander estimator of a monotone density are stated by Anevski and Hössjer (2006). We find that the respective limiting distributions of the MLE and its first derivative depend on a stochastic process, $H_k$, equal almost surely to the "lower invelope" (or just "invelope") on $\mathbb{R}$ of the integrated Brownian motion minus $t^{k+2}$, where $k$ is the order of the first nonzero derivative of $\varphi_0$ at the point of interest.

In Section 3, the estimation point $x_0$ is taken to be equal to the mode, $m_0$, defined to be the smallest point in the modal interval of the log-concave density $f_0$. A natural estimator of $m_0$, which we denote by $\widehat{M}_n$, can be taken to be the smallest number maximizing the MLE $\widehat{\varphi}_n$ or, equivalently, the smallest number maximizing the MLE $\widehat{f}_n$. In this section, we establish our second main result: the asymptotic distribution of $\widehat{M}_n$. Under the assumption that the second derivative $f_0''(m_0) < 0$, we show that this distribution depends on the random variable defined to be the argmax or mode of $H_2^{(2)}$ on $\mathbb{R}$. When the second, third and higher derivatives of order $k-1$ or lower vanish at $m_0$ but $f_0^{(k)}(m_0) < 0$, then the limit distribution depends on the mode of $H_k^{(2)}$.

Proofs are deferred to Section 4.

To illustrate all the quantities for which we provide limiting distributions, in Figure 1 we give plots of $\widehat{f}_n$, $\widehat{\varphi}_n$, $\widehat{F}_n$ and $\widehat{\lambda}_n = \widehat{f}_n/(1 - \widehat{F}_n)$, based on two samples of sizes $n = 20$ and $n = 200$ drawn from a Gamma$(2,1)$ density $f_0(x) = xe^{-x}1_{[0,\infty)}(x)$. All these plots were generated using the R-package logcondens [see Rufibach and Dümbgen (2007)].

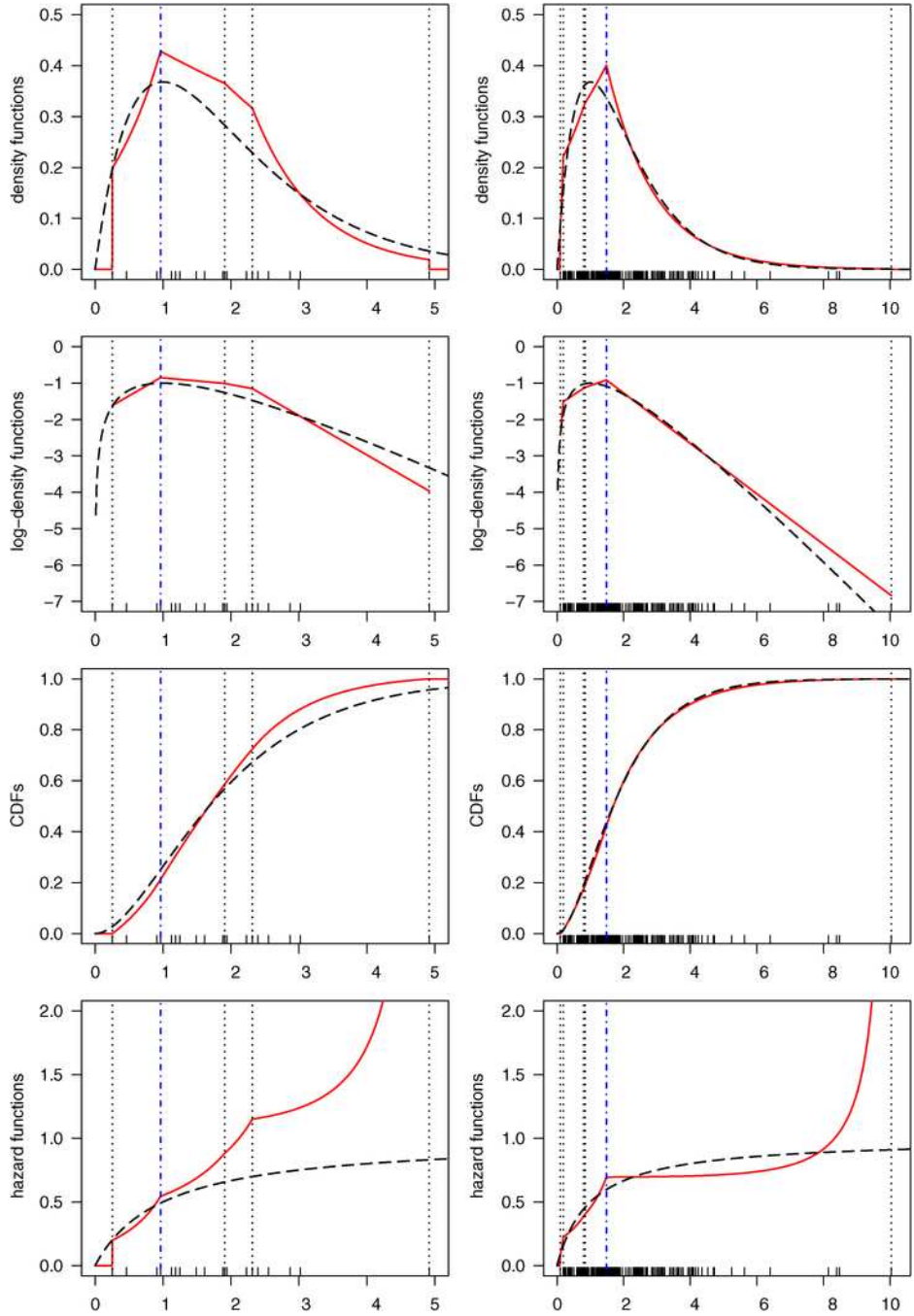

FIG. 1. *Examples for log-concave density, log-density, CDF, and hazard rate estimation for $n = 20, 200$ ($--$ true functions, $-$ estimators). The dotted vertical lines indicate the set $\widehat{\mathcal{S}}_n(\widehat{\varphi}_n)$. The $\cdot - \cdot -$ vertical lines are placed at the mode of the estimated density.*

**2. Limiting distribution theory.** To state the main result, we make the following assumptions.

2.1. *Assumptions.* Fix $x_0 \in \mathbb{R}$. We suppose that the true density $f_0 = \exp \varphi_0$ satisfies the following assumptions:

(A1) The density function $f_0 \in \mathcal{LC}$.
(A2) $f_0(x_0) > 0$.
(A3) The function $\varphi_0$ is at least twice continuously differentiable in a neighborhood of $x_0$.
(A4) If $\varphi_0''(x_0) \neq 0$, then $k = 2$. Otherwise, suppose that $k$ is the smallest integer such that $\varphi_0^{(j)}(x_0) = 0, j = 2, \ldots, k-1$, and $\varphi_0^{(k)}(x_0) \neq 0$, and $\varphi_0^{(k)}$ is continuous in a neighborhood of $x_0$.

Note that concavity of $\varphi_0$ and (A3) and (A4) imply that $k$ is necessarily even and that $\varphi_0^{(k)}(x_0) < 0$. Indeed, suppose that $k > 2$. Using Taylor expansion of $\varphi_0''$ up to degree $k-2$, there exists a small $h > 0$ for which we can write

$$\varphi_0''(x) = \frac{\varphi_0^{(k)}(x_0)}{(k-2)!}(x-x_0)^{k-2} + o((x-x_0)^{k-2}), \qquad x \in [x_0 - h, x_0 + h].$$

Since $\varphi_0''(x) \leq 0$ for all $x \in [x_0 - h, x_0 + h]$, it follows that $k-2$ is even [i.e., $k$ is even and $\varphi_0^{(k)}(x_0) < 0$].

2.2. *Notation.* Let $W$ denote two-sided Brownian motion, starting at 0. For $t \in \mathbb{R}$, define:

$$(2.1) \qquad Y_k(t) = \begin{cases} \int_0^t W(s)\,ds - t^{k+2}, & \text{if } t \geq 0, \\ \int_t^0 W(s)\,ds - t^{k+2}, & \text{if } t < 0. \end{cases}$$

For the uniform norm of a bounded function $f$, we write $\|f\|_\infty = \sup_{x \in \mathbb{R}} |f(x)|$. The derivative of $\widehat{\varphi}_n$ at $x \in \mathbb{R}$ is as usual denoted by $\widehat{\varphi}_n'(x)$. However, if $x \in \widehat{\mathcal{S}}_n(\widehat{\varphi}_n)$, then we define $\widehat{\varphi}_n'(x)$ as the left-derivative.

THEOREM 2.1. *Suppose that* (A1)–(A4) *hold. Then,*

$$\begin{pmatrix} n^{k/(2k+1)}(\widehat{f}_n(x_0) - f_0(x_0)) \\ n^{(k-1)/(2k+1)}(\widehat{f}_n'(x_0) - f_0'(x_0)) \end{pmatrix} \xrightarrow{d} \begin{pmatrix} c_k(x_0, \varphi_0) H_k^{(2)}(0) \\ d_k(x_0, \varphi_0) H_k^{(3)}(0) \end{pmatrix}$$

*and*

$$\begin{pmatrix} n^{k/(2k+1)}(\widehat{\varphi}_n(x_0) - \varphi_0(x_0)) \\ n^{(k-1)/(2k+1)}(\widehat{\varphi}_n'(x_0) - \varphi_0'(x_0)) \end{pmatrix} \xrightarrow{d} \begin{pmatrix} C_k(x_0, \varphi_0) H_k^{(2)}(0) \\ D_k(x_0, \varphi_0) H_k^{(3)}(0) \end{pmatrix},$$

*where $H_k$ is the "lower invelope" of the process $Y_k$; that is,*

$H_k(t) \leq Y_k(t)$ for all $t \in \mathbb{R}$;

$H_k^{(2)}$ is concave;

$H_k(t) = Y_k(t)$, if the slope of $H_k^{(2)}$ decreases strictly at $t$.

The constants $c_k$, $d_k$, $C_k$ and $D_k$ are given by

$$c_k(x_0, \varphi_0) = \left( \frac{f_0(x_0)^{k+1} |\varphi_0^{(k)}(x_0)|}{(k+2)!} \right)^{1/(2k+1)}, \tag{2.2}$$

$$d_k(x_0, \varphi_0) = \left( \frac{f_0(x_0)^{k+2} |\varphi_0^{(k)}(x_0)|^3}{[(k+2)!]^3} \right)^{1/(2k+1)}, \tag{2.3}$$

$$C_k(x_0, \varphi_0) = \left( \frac{|\varphi_0^{(k)}(x_0)|}{f_0(x_0)^k (k+2)!} \right)^{1/(2k+1)}, \tag{2.4}$$

$$D_k(x_0, \varphi_0) = \left( \frac{|\varphi_0^{(k)}(x_0)|^3}{f_0(x_0)^{k-1} [(k+2)!]^3} \right)^{1/(2k+1)}. \tag{2.5}$$

COROLLARY 2.2. *Suppose that* (A1)–(A4) *hold with* $k=2$. *Then,*

$$\begin{pmatrix} n^{2/5}(\widehat{f}_n(x_0) - f_0(x_0)) \\ n^{1/5}(\widehat{f}'_n(x_0) - f'_0(x_0)) \end{pmatrix} \xrightarrow{d} \begin{pmatrix} c_2(x_0, \varphi_0) H_2^{(2)}(0) \\ d_2(x_0, \varphi_0) H_2^{(3)}(0) \end{pmatrix}$$

*and*

$$\begin{pmatrix} n^{2/5}(\widehat{\varphi}_n(x_0) - \varphi_0(x_0)) \\ n^{1/5}(\widehat{\varphi}'_n(x_0) - \varphi'_0(x_0)) \end{pmatrix} \xrightarrow{d} \begin{pmatrix} C_2(x_0, \varphi_0) H_2^{(2)}(0) \\ D_2(x_0, \varphi_0) H_2^{(3)}(0) \end{pmatrix},$$

*where* $H_2$ *is the (concave) invelope of the process* $Y_2$; *that is,*

$H_2(t) \leq Y_2(t)$ for all $t \in \mathbb{R}$;

$H_2^{(2)}$ is concave;

$H_2(t) = Y_2(t)$ if the slope of $H_2^{(2)}$ decreases strictly at $t$.

*The constants* $c_2$, $d_2$, $C_2$ *and* $D_2$ *are given by* (2.2)–(2.5), *with* $k=2$.

Note that the constants $C_2(x_0, \varphi_0)$ and $D_2(x_0, \varphi_0)$, up to inversion of $f_0(x_0)$, exhibit a structure very similar to that of the constants given by Groeneboom, Jongbloed and Wellner (2001b) in the problem of estimating a convex density $g_0$ on $[0, \infty)$. We recall here that, in the latter problem, those constants are found to be equal to (we use our notation to make the comparison easy)

$$c_2(x_0, g_0) = \left( \frac{g_0(x_0)^2 g_0^{(2)}(x_0)}{4!} \right)^{1/5}, \qquad d_2(x_0, g_0) = \left( \frac{g_0(x_0)(g_0^{(2)}(x_0))^3}{(4!)^3} \right)^{1/5}.$$

It is clear that $\varphi_0$ in the log-concave problem plays exactly the same role as $f_0$ in the problem of estimating a convex density. However, in the first case estimation is based on observations which are distributed according to $\exp \varphi_0$, whereas in the latter the data come from $f_0$ itself. A good insight into the difference between the expressions of the asymptotic constants can be gained from the proof of Theorem 4.6 in Section 4. There, we show that the leading coefficient of the drift of the limiting process $Y_k$ depends on $\varphi_0^{(k)}(x_0) f_0(x_0) = f_0^{(k)}(x_0) - (\varphi_0'(x_0))^k f_0(x_0)$, where the second term is "filtered out" in the Taylor expansion of the estimation error in the neighborhood of $x_0$. Hence, $|\varphi_0^{(k)}(x_0)| \cdot f_0(x_0)$ can be viewed as the dominating term replacing $|g_0^{(k)}(x_0)|$ in the convex estimation problem. For $k=2$, the constants $c_2(x_0, \varphi_0)$ and $d_2(x_0, \varphi_0)$ given in (2.2) and (2.3), with $k=2$, match closely with $c_2(x_0, g_0)$ and $d_2(x_0, g_0)$ obtained by Groeneboom, Jongbloed and Wellner (2001b) in the convex estimation problem, with $f_0(x_0)$ in the numerator, whereas $f_0(x_0)$ shows up in the denominator in the asymptotic constants $C_2(x_0, \varphi_0)$ and $D_2(x_0, \varphi_0)$. This results from applying the delta-method to $\widehat{f}_n(x_0) = \exp(\widehat{\varphi}_n(x_0))$ and $\widehat{f}'_n(x_0) = \widehat{\varphi}'_n(x_0) \widehat{f}_n(x_0)$, which yields $C_2(x_0, \varphi_0)$ and $D_2(x_0, \varphi_0)$.

Here is an explicit example showing how vanishing second (and higher) derivatives can occur. Consider the density function

$$f_0(x) = \sqrt{2} \frac{\Gamma(3/4)}{\pi} \exp(-x^4), \qquad x \in \mathbb{R}.$$

In this case $\varphi_0^{(j)}(x_0) = 0, j = 1, 2, 3$ for $x_0 = 0$, and $\varphi_0^{(4)}(x_0) \neq 0$. The following "tilted" version of $f_0$ shows that vanishing second derivatives of $\varphi_0$ can also occur at points other than the mode of $f$:

$$\tilde{f}_0(x) = \exp(a + bx) f_0(x) = \tilde{a} \exp(bx - x^4),$$

where $\tilde{a} = \tilde{a}(b) := 1/\int_{\mathbb{R}} \exp(bx - x^4)\, dx$; in this case, $\tilde{\varphi}_0 := \log \tilde{f}_0$ satisfies $\tilde{\varphi}_0''(0) = 0$, but the mode $\tilde{m}_0 := M(\tilde{f}_0) = (b/4)^{1/3} > 0$ when $b > 0$, and $\tilde{\varphi}_0''(\tilde{m}_0) = -12(b/4)^{2/3} < 0$.

Finally, and in order to compare also the random parts of the limits in the convex and log-concave estimation problems, we would like to note that for our lower invelope process $H_k$, $-H_k$ has the same distribution as the "upper invelope" of $-Y_k$, which was called just the "invelope" in the case $k=2$ by Groeneboom, Jongbloed and Wellner (2001b): The process $-Y_k$ has a drift equal to plus $t^{k+2}$, which specializes to $t^4$ in the convex density problem with $k=2$. This "upper invelope" stays above $-Y_k$ and admits a convex second derivative. Since $-W$ has the same distribution as $W$, it follows that the upper and lower invelopes $\mathbb{H}_k$ and $H_k$ (associated with estimation of convex and concave functions, resp.) satisfy $\mathbb{H}_k \stackrel{d}{=} -H_k$. Since the derivatives at zero $\mathbb{H}_k^{(2)}(0)$ and $\mathbb{H}_k^{(3)}(0)$ of $\mathbb{H}_k$ are distributed symmetrically about zero, the same is true of the derivatives at zero $H_k^{(2)}(0)$ and $H_k^{(3)}(0)$ of $H_k$.

As shown by Barlow and Proschan (1975), Lemma 5.8, page 77 [see also Marshall and Olkin (1979), page 493; Marshall and Olkin (2007), page 102; An (1998) and Bagnoli and Bergstrom (2005)], if $f_0$ is log-concave, then the hazard function

$$\lambda_0(x) = \frac{f_0(x)}{1 - F_0(x)} 1_{\{x < F_0^{-1}(1)\}}$$

is monotone nondecreasing. Defining the estimator of $\lambda_0$ based on $\widehat{f}_n$ as

$$\widehat{\lambda}_n(x) = \frac{\widehat{f}_n(x)}{1 - \widehat{F}_n(x)} 1_{\{x < X_{(n)}\}},$$

application of the delta-method yields the following corollary.

COROLLARY 2.3. *Suppose that* (A1)–(A4) *hold. Then,*

$$\begin{pmatrix} n^{k/(2k+1)}(\widehat{\lambda}_n(x_0) - \lambda_0(x_0)) \\ n^{(k-1)/(2k+1)}(\widehat{\lambda}'_n(x_0) - \lambda'_0(x_0)) \end{pmatrix} \xrightarrow{d} \begin{pmatrix} g_k(x_0, \varphi_0) H_k^{(2)}(0) \\ h_k(x_0, \varphi_0) H_k^{(3)}(0) \end{pmatrix},$$

*where the constants $g_k$ and $h_k$ are given by*

$$g_k(x_0, \varphi_0) = c_k(x_0, \varphi_0)/(1 - F_0(x_0))$$
$$h_k(x_0, \varphi_0) = d_k(x_0, \varphi_0)/(1 - F_0(x_0)).$$

**3. Inference about the mode of $f_0$.** Estimation of the mode of a unimodal density has been considered by many authors [see, e.g., Parzen (1962), Chernoff (1964), Grenander (1965), Dalenius (1965), Venter (1967), Wegman (1970a, 1970b, 1971), Eddy (1980, 1982), Hall (1982), Müller (1989), Romano (1988), Vieu (1996) and, more recently, Meyer (2001) and Herrmann and Ziegler (2004)].

Empirical studies of the performance of various estimators are given by Dalenius (1965), Ekblom (1972), Meyer (2001) and Meyer and Woodroofe (2004). Many of the methods considered for estimating the mode of a unimodal smooth density use kernel estimation, but others are based on the principle of substitution with another choice of estimator of the population density. For example, the estimators of Venter (1967) are related to nearest-neighbor estimators of the density $f_0$. All the estimators of the mode in the class of unimodal densities known to us involve some more or less ad hoc choice, essentially because the maximum likelihood estimator of a unimodal density is not well defined, as explained by Birgé (1997). [Note that Wegman (1970b, 1971) discussed the nonparametric MLE of a unimodal density subject to a constraint on the height of the mode; without some constraint of this type, the MLE does not exist.]

For virtually all of the estimators of which we are aware, some choice of a smoothing parameter, bandwidth or constraint is required. Empirical choice of smoothing parameters has been studied by Müller (1989), who studied local methods of choosing the smoothing parameter, Grund and Hall (1995), who studied bootstrap methods, and Ziegler (2004), who studied plug-in methods. Klemelä (2005) gave a construction of adaptive estimators based on Lepski's method [Lepskiĭ (1991, 1992)]. For nonparametric Bayes estimators of unimodal densities and, hence, of the mode [see Brunner and Lo (1989) and Ho (2006a, 2006b)]; for these estimators, choice of a prior is equivalent to a choice of smoothing parameters.

In contrast, estimation in the (large) subclass of log-concave (or strongly unimodal) densities is much simpler, avoiding bandwidth or smoothing parameter choices completely. Since the maximum likelihood estimator exists, we can simply estimate the mode by the mode (or smallest point in a modal interval) of the MLE $\widehat{f}_n$. Using the notation introduced by Eddy (1982) [and also used by Romano (1988)], we let $\widehat{M}_n := M(\widehat{f}_n)$ where $M$ denotes the mode functional (or "smallest argmax" functional) given by

$$M(g) := \min\left\{t : g(t) = \max_{u \in \mathbb{R}} g(u)\right\}.$$

Because of the adaptive properties of the MLE's $\widehat{f}_n$ of $f_0$ and $\widehat{\varphi}_n$ of $\varphi_0$, discussed in Section 1, we expect $\widehat{M}_n$ to adapt to different local smoothness (or peakedness) hypotheses on $f_0$ [much as the Grenander estimator is locally adaptive in the case of estimating a monotone density, see, for example, Birgé (1989), page 1535]. Here, we study $\widehat{M}_n$ as an estimator of the mode $M(f_0) := m_0$ under just the condition that $f_0$ has a continuous second derivative $f_0''$ in a neighborhood of $m_0$, with $f_0''(m_0) < 0$. We begin in the next subsection with a new asymptotic minimax lower bound for estimation of $m_0$ under this hypothesis. The following subsection gives our new limiting distribution result for the MLE $\widehat{M}_n$ of the mode $m_0$.

3.1. *New lower bounds for estimating the mode.* Has'minskiĭ (1979) established a lower bound for estimation of the mode $m_0$ of a unimodal density $f \in \mathcal{U}$, assuming that $f$ satisfies $f''(m_0) < 0$. He showed that the best local asymptotic minimax rate of convergence for *any* estimator of $m_0$ is $n^{-1/5}$. Has'minskiĭ based his proof on a sequence of parametric submodels of the form

$$f_n(x, \theta) = f(x) + \theta n^{-2/5} g(n^{1/5}(x - m_0)),$$

where, for $a := -f''(m_0)$,

$$g(x) := g_a(x) = \begin{cases} x, & \text{if } |x| \leq 1/a, \\ 0, & \text{if } |x| \geq K > 1/a \end{cases}$$

and $g := g_a$ satisfies $g(-x) = -g(x)$ and $|g''(x)| < a/2$ for all $x \in \mathbb{R}$. However, Has'minskiĭ (1979) did not study the dependence of the local minimax bound on $a = -f''(m_0)$ and $f(m_0)$, leaving his bound in terms of $c_0^2 := f(m_0)/\int g_a^2(x)\,dx$ involving the still unspecified function $g = g_a$.

Here, we consider different parametric submodels and derive the dependence of the constant in local asymptotic minimax lower bound for estimation of the mode $m_0$ in the family $\mathcal{LC}$ of log-concave (or strongly unimodal) densities.

We want to derive asymptotic lower bounds for the local minimax risks for estimating the mode $M(f)$. The $L_1$-minimax risk for estimating a functional $\nu$ of $f_0$, based on a sample $X_1, \ldots, X_n$ of size $n$ from $f_0$, which is known to be in a subset $\mathcal{LC}_{n,\tau}$ of $\mathcal{LC}$ is defined by

$$(3.1) \qquad MMR_1(n, T_n, \mathcal{LC}_{n,\tau}) := \inf_{T_n} \sup_{f \in \mathcal{LC}_{n,\tau}} E_f|T_n - \nu(f)|,$$

where the infimum ranges over all possible measurable functions $T_n = t_n(X_1, \ldots, X_n)$ mapping $\mathbb{R}^n$ to $\mathbb{R}$. The shrinking classes $\mathcal{LC}_{n,\tau}$ used here are Hellinger balls centered at $f_0$:

$$\mathcal{LC}_{n,\tau} = \left\{ f \in \mathcal{LC} : H^2(f, f_0) = \tfrac{1}{2} \int_{-\infty}^{\infty} (\sqrt{f(z)} - \sqrt{f_0(z)})^2 \, dz \leq \tau/n \right\}.$$

Consider estimation of

$$(3.2) \qquad \nu(f) := M(f) = \inf\left\{ t \in \mathbb{R} : t = \sup_{u \in \mathbb{R}} f(u) \right\}.$$

Let $f_0 \in \mathcal{LC}$ and $m_0 = M(f_0)$ be fixed, such that $f_0$ is twice continuously differentiable at $m_0$ and $f_0''(m_0) < 0$. Consider the family $\{\varphi_\varepsilon\}_{\varepsilon > 0}$ and resulting family $\{f_\varepsilon\}_{\varepsilon > 0}$, defined as follows much as:

$$\varphi_\varepsilon(x) = \begin{cases} \varphi_0(x), & x < m_0 - \varepsilon c_\varepsilon, \\ \varphi_0(x), & x > m_0 + \varepsilon, \\ \varphi_0(m_0 + \varepsilon), \\ \quad + \varphi_0'(m_0 + \varepsilon)(x - m_0 - \varepsilon), & x \in [m_0 - \varepsilon, m_0 + \varepsilon], \\ \varphi_0(m_0 - \varepsilon c_\varepsilon), \\ \quad + \varphi_0'(m_0 - \varepsilon c_\varepsilon)(x - m_0 + \varepsilon c_\varepsilon), & x \in [m_0 - \varepsilon c_\varepsilon, m_0 - \varepsilon), \end{cases}$$

where $c_\varepsilon$ is chosen so that $\varphi_\varepsilon$ is continuous at $m_0 - \varepsilon$. Note that if $\varphi_0(x) = \gamma - \gamma_0(x - m_0)^2$, then $c_\varepsilon = 3$, for all $\varepsilon$, and $c_\varepsilon \to 3$, as $\varepsilon \downarrow 0$, since $f_0''(m_0) < 0$. Now define

$$h_\varepsilon(x) := \exp(\varphi_\varepsilon(x)) \quad \text{and} \quad f_\varepsilon(x) := \frac{h_\varepsilon(x)}{\int h_\varepsilon(y)\,dy}.$$

Then, $f_\varepsilon$ is log-concave for each $\varepsilon > 0$ with mode $m_0 - \varepsilon$ by construction, so with $\nu(f_\varepsilon) := M(f_\varepsilon) :=$ the mode of $f_\varepsilon$, we have

$$\nu(f_\varepsilon) - \nu(f_0) = M(f_\varepsilon) - M(f_0) = m_0 - \varepsilon - m_0 = -\varepsilon.$$

Furthermore, the following lemma holds.

LEMMA 3.1. *Under the above assumptions,*

$$H^2(f_\varepsilon, f_0) = \frac{2 f_0''(m_0)^2}{5 f_0(m_0)} \varepsilon^5 + o(\varepsilon^5) := \rho \varepsilon^5 + o(\varepsilon^5).$$

PROOF. Proceeding as in Jongbloed (1995),

$$H^2(f_\varepsilon, f_0) = \frac{1}{2} \int_{-\infty}^{\infty} [\sqrt{f_\varepsilon(x)} - \sqrt{f_0(x)}]^2 \, dx$$

$$= \frac{1}{2} \int_{m_0 - \varepsilon c_\varepsilon}^{m_0 + \varepsilon} [\sqrt{f_\varepsilon(x)} - \sqrt{f_0(x)}]^2 \, dx$$

$$= \frac{2}{5} f_0(m_0) \varphi_0''(m_0)^2 \varepsilon^5 + o(\varepsilon^5) = \frac{2}{5} \frac{f_0''(m_0)^2}{f_0(m_0)} \varepsilon^5 + o(\varepsilon^5)$$

as $\varepsilon \downarrow 0$. Calculations similar to those of Jongbloed (1995) [see also Jongbloed (2000) and Groeneboom, Jongbloed and Wellner (2001b)] complete the proof of the lemma. □

Taking $\varepsilon = c n^{-1/5}$ and defining $f_n := f_{cn^{-1/5}}$ yields

$$\nu(f_n) - \nu(f_0) = M(f_n) - M(f_0) = -c n^{-1/5}$$

and

$$n H^2(f_n, f_0) = \frac{2}{5} \frac{f_0''(m_0)^2}{f_0(m_0)} c^5 + o(1) := \rho c^5 + o(1).$$

Plugging these into the lower bound Lemma 4.1 of Groeneboom (1996), with $\ell(x) := |x|$, yields

$$\liminf_n \inf_{T_n} n^{1/5} \max\{E_{n,P_n}|T_n - M(f_n)|, E_{n,P}|T_n - M(f_0)|\}$$

$$\geq \frac{1}{4} c \exp(-2\rho c^5) = \frac{e^{-1/5}}{4 \cdot 10^{1/5}} \rho^{-1/5} = (0.15512) \left( \frac{f_0(m_0)}{f_0''(m_0)^2} \right)^{1/5}$$

by choosing $c = (10\rho)^{-1/5}$. This yields the following proposition.

PROPOSITION 3.2 (Minimax risk lower bound). *Suppose that $\nu(f) = M(f)$, as defined in (3.2), and that $\mathcal{LC}_{n,\tau}$ is as defined above where $f_0''$ is continuous in a neighborhood of $m_0 = M(f_0)$ with $f_0''(m_0) < 0$. Then,*

$$\sup_{\tau > 0} \limsup_{n \to \infty} n^{1/5} \inf_{T_n} \sup_{f \in \mathcal{LC}_{n,\tau}} E_f |T_n - M(f)|$$

$$\geq \left(\frac{5/2}{4^5 \cdot e \cdot 10}\right)^{1/5} \left(\frac{f_0(m_0)}{f_0''(m_0)^2}\right)^{1/5} \doteq (0.15512) \left(\frac{f_0(m_0)}{f_0''(m_0)^2}\right)^{1/5}.$$

REMARK 3.3. Note that the constant $b(f_0, m_0) := (f_0(m_0)/f_0''(m_0)^2)^{1/5}$ appearing on the right-hand side of this lower bound is scale equivariant in exactly the right way: if $f_c(x) := f_0(m_0 + (x - m_0)/c)/c$ for $c > 0$, then $b(f_c, m_0) = cb(f_0, m_0)$ for all $c > 0$. The constant $b(f_0, m_0)$ will appear in the limit distribution appearing in the next subsection.

REMARK 3.4. If $\mathcal{LC}$ is replaced by the class $\mathcal{U}$ of unimodal densities on $\mathbb{R}$ and $\mathcal{LC}_{n,\tau}$ is replaced by $\mathcal{U}_{n,\tau}$ defined analogously where $f_0$ satisfies $f_0''(m_0) < 0$ and $f_0''$ continuous in a neighborhood of $m_0$, then a minimax lower bound of the same form as Proposition 3.2 holds with exactly the same dependence on $b(f_0, m_0) = (f_0(m_0)/f_0''(m_0)^2)^{1/5}$, but with the absolute constant $0.15512\ldots$ replaced by $0.19784\ldots$. This can be seen by taking the perturbations $\{f_\varepsilon\}_{\varepsilon > 0}$ defined by

$$f_\varepsilon(x) = \begin{cases} f_0(x), & x \leq x_0 - \varepsilon, \\ f_0(x), & x > x_0 + \varepsilon, \\ f_0(x_0) + b_\varepsilon(x - x_0 + \varepsilon), & x_0 - \varepsilon \leq x \leq x_0 + \varepsilon, \end{cases}$$

where $b_\varepsilon$ is chosen so that $f_\varepsilon(x_0 + \varepsilon) > f_0(x_0 + \varepsilon)$ and $\int_{x_0 - \varepsilon}^{x_0 + \varepsilon} f_\varepsilon(x)\,dx = \int_{x_0 - \varepsilon}^{x_0 + \varepsilon} f_0(x)\,dx$.

REMARK 3.5. If $\varphi_0$ is continuously $k$-times differentiable in a neighborhood of the mode $m_0$, $\varphi_0^{(j)}(m_0) = 0$ for $j = 2, \ldots, k - 1$ and $\varphi_0^{(k)}(m_0) \neq 0$ [assumption (A4)], then it can be shown that the minimax rate of convergence is $n^{1/(2k+1)}$ and that the minimax lower bound is proportional to

$$\left(\frac{1}{f_0(m_0)\varphi_0^{(k)}(m_0)^2}\right)^{1/(2k+1)} = \left(\frac{f_0(m_0)}{f_0^{(k)}(m_0)^2}\right)^{1/(2k+1)},$$

where the proportionality constant depends on the largest root of the polynomial $x^k - (k/(k-1))x^{k-1} - (2k-1)/(k-1)$ (which equals 3 when $k = 2$).

3.2. *Limiting distribution for the MLE $\widehat{M}_n$ in $\mathcal{LC}$.* Now, let $\widehat{f}_n$ be the MLE of $f$ in the class $\mathcal{LC}$ of log-concave densities, and let $\widehat{M}_n = M(\widehat{f}_n)$, $m_0 = M(f_0)$. Here is our result concerning the limiting distribution of $\widehat{M}_n$ under the same assumptions on $f_0$ as in the previous section on lower bounds.

THEOREM 3.6. *Suppose that $f_0''$ is continuous in a neighborhood of $m_0 = M(f_0)$ and that $f_0''(m_0) < 0$. Then,*

$$n^{1/5}(\widehat{M}_n - m_0) \xrightarrow{d} \left(\frac{(4!)^2 f_0(m_0)}{f_0''(m_0)^2}\right)^{1/5} M(H_2^{(2)}).$$

Note that the limiting distribution depends on a multiple of the same constant $b(f_0, m_0)$, which appears in the asymptotic minimax lower bound of Proposition 3.2, times a universal term $M(H_2^{(2)})$, the mode of the "estimator" $H_2^{(2)}(t)$ of the canonical concave function $-12t^2$ in the limit Gaussian problem: estimate the mode of $f_0(t) = -12t^2$, based on observation of $Y(t) = \int_0^t X(s)\, ds$, when

$$dX(t) = f_0(t)\, dt + dW(t).$$

We expect that this distribution, namely the distribution of

$$M(H_2^{(2)}) = \arg\max_{t \in \mathbb{R}} H_2^{(2)}(t),$$

will occur in several other problems involving nonparametric estimation of the mode or antimode of convex or concave functions under similar second derivative hypotheses. For example, it seems clear that it will occur as the limiting distribution of the nonparametric estimator of the antimode of a convex bathtub-shaped hazard [in the setting of Jankowski and Wellner (2007)]; as the limiting distribution of the nonparametric estimator of the antimode of a convex regression function in the setting of Groeneboom, Jongbloed and Wellner (2001b); and as the limiting distribution of the nonparametric estimator of the mode of a concave regression function.

When $\varphi_0^{(j)}(m_0) = 0$, for $j = 2, \ldots, k-1$, $\varphi_0^{(k)}(m_0) \neq 0$, and $\varphi_0^{(k)}$ is continuous in a neighborhood of $m_0$, then an analogous result (with a completely similar proof) holds:

$$n^{1/(2k+1)}(\widehat{M}_n - m_0) \xrightarrow{d} \left(\frac{(k+2)!^2}{f_0(m_0)|\varphi_0^{(k)}(m_0)|^2}\right)^{1/(2k+1)} M(H_k^{(2)}).$$

In particular, when $k = 4$, the rate of convergence is $n^{1/9}$, and the limit distribution becomes that of

$$\left(\frac{6!^2 f_0(m_0)}{f_0^{(4)}(m_0)^2}\right)^{1/9} M(H_4^{(2)}).$$

Apparently, estimation of $m_0$ becomes considerably more difficult when the second and possibly higher order derivatives of $\varphi_0$ vanish at $m_0$.

On the other hand, if $\varphi_0$ (or equivalently, $f_0$) is cusp-shaped at $m_0$, then the rate of convergence of $\widehat{M}_n$ is $n^{1/3}$, and the local asymptotic minimax rate of convergence is also $n^{1/3}$; we will pursue these issues elsewhere.

## 4. Proofs for Sections 2 and 3.

Throughout this section, we fix $k$ and let

$$r_n := n^{(k+2)/(2k+1)}, \qquad s_n := n^{-1/(2k+1)},$$

$$x_n(t) := x_{n,k}(t) := x_0 + s_n t := x_0 + n^{-1/(2k+1)} t,$$

$$\mathcal{I} := \mathcal{I}(x_0, n, k, t) := \begin{cases} [x_0, x_n(t)], & t \geq 0, \\ [x_n(t), x_0], & t < 0. \end{cases}$$

4.1. *Preparation: technical lemmas and tightness results.* First, some notation.

*Local processes*: The local processes $\mathbb{Y}_n^{\mathrm{loc}}$ and $\widehat{H}_n^{\mathrm{loc}}$ are defined for $t \in \mathbb{R}$ by

$$\mathbb{Y}_n^{\mathrm{loc}}(t) := r_n \int_{x_0}^{x_n(t)} \left( \mathbb{F}_n(v) - \mathbb{F}_n(x_0) - \int_{x_0}^{v} \left( \sum_{j=0}^{k-1} \frac{f_0^{(j)}(x_0)}{j!} (u - x_0)^j \right) du \right) dv$$

and

$$\widehat{H}_n^{\mathrm{loc}}(t) := r_n \int_{x_0}^{x_n(t)} \int_{x_0}^{v} \left( \widehat{f}_n(u) - \sum_{j=0}^{k-1} \frac{f_0^{(j)}(x_0)}{j!} (u - x_0)^j \right) du\, dv$$

$$+ \hat{A}_n t + \hat{B}_n,$$

where in the limit Gaussian problem: estimate the mode

(4.1) $$\hat{A}_n = r_n s_n (\widehat{F}_n(x_0) - \mathbb{F}_n(x_0)) \quad \text{and}$$

(4.2) $$\hat{B}_n = r_n (\widehat{H}_n(x_0) - \mathbb{H}_n(x_0)).$$

We also define the "modified" local processes

(4.3)
$$\mathbb{Y}_n^{\mathrm{locmod}}(t) := \frac{r_n}{f_0(x_0)} \int_{x_0}^{x_n(t)} \Bigg( \mathbb{F}_n(v) - \mathbb{F}_n(x_0)$$
$$- \int_{x_0}^{v} \left( \sum_{j=0}^{k-1} \frac{f_0^{(j)}(x_0)}{j!} (u - x_0)^j \right) du \Bigg) dv$$
$$- r_n \int_{x_0}^{x_n(t)} \int_{x_0}^{v} \widehat{\Psi}_{k,n,2}(u)\, du\, dv$$

and

$$\mathbb{H}_n^{\text{locmod}}(t) := r_n \int_{x_0}^{x_n(t)} \int_{x_0}^{v} (\widehat{\varphi}_n(u) - \varphi_0(x_0) - (u - x_0)\varphi_0'(x_0))\, du\, dv$$
(4.4)
$$+ \frac{\hat{A}_n t + \hat{B}_n}{f_0(x_0)},$$

where $\widehat{\Psi}_{k,n,2}$ is defined below in (4.26).

The following lemma uses the notion of uniform covering numbers [see van der Vaart and Wellner (1996), Sections 2.1 and 2.7] for complete definitions and further information.

LEMMA 4.1. *Let $\mathcal{F}$ be a collection of functions defined on $[x_0 - \delta, x_0 + \delta]$, with $\delta > 0$ small and let $s > 0$. Suppose that for a fixed $x \in [x_0 - \delta, x_0 + \delta]$ and $R > 0$, such that $[x, x+R] \subseteq [x_0 - \delta, x_0 + \delta]$, the collection*

$$\mathcal{F}_{x,R} = \{f_{x,y} := f 1_{[x,y]}, f \in \mathcal{F}, x \leq y \leq x + R\}$$

*admits an envelope $F_{x,R}$, such that*

$$EF_{x,R}^2(X_1) \leq KR^{2d-1}, \qquad R \leq R_0$$

*for some $d \geq 1/2$ and $K > 0$, depending only on $x_0$ and $\delta$. Moreover, suppose that*

(4.5) $$\sup_Q \int_0^1 \sqrt{\log N(\eta \|F_{x,R}\|_{Q,2}, \mathcal{F}_{x,R}, L_2(Q))}\, d\eta < \infty.$$

*Then, for each $\varepsilon > 0$, there exist random variables $M_n$ of order $O_p(1)$ (not depending on $x$ or $y$) and $R_0 > 0$, such that*

$$\left| \int f_{x,y}\, d(\mathbb{F}_n - F_0) \right| \leq \varepsilon |y - x|^{s+d} + n^{-(s+d)/(2s+1)} M_n \qquad \text{for } |y - x| \leq R_0.$$

PROOF. See Kim and Pollard (1990) and Balabdaoui and Wellner (2007), Lemmas 4.4 and 6.1. The special case $s = 1 = d$ is Lemma 4.1 of Kim and Pollard (1990). □

LEMMA 4.2. *If (A3) and (A4) hold, then*

(4.6) $$f_0^{(j)}(x_0) = [\varphi_0'(x_0)]^j f_0(x_0) \qquad \text{for } j = 1, \ldots, k-1$$

*and, for $j = k$*

$$f_0^{(k)}(x_0) = (\varphi_0^{(k)}(x_0) + [\varphi_0'(x_0)]^k) f_0(x_0).$$

PROOF. The expressions for $f_0^{(j)}(x_0)$ follow immediately from a recursive argument using the identity $f_0 = \exp \varphi_0$ and the assumption $\varphi_0^{(j)}(x_0) = 0$, for $j = 2, \ldots, k-1$, if $k > 2$. □

Now, let $\tau_n^+ := \inf\{t \in \widehat{\mathcal{S}}(\widehat{\varphi}_n) : t > x_0\}$ and $\tau_n^- := \sup\{t \in \widehat{\mathcal{S}}(\widehat{\varphi}_n) : t < x_0\}$.

THEOREM 4.3. *If* (A1)–(A4) *hold, then*

$$\tau_n^+ - \tau_n^- = O_p(n^{-1/(2k+1)}). \tag{4.7}$$

Theorem 4.3 should be compared to Theorem 3.3 of Dümbgen and Rufibach (2009). When their Theorem 3.3 is specialized to the case $\beta = 2$, so that $\varphi_0''(x) \leq C < 0$, for all $x \in T := [A, B]$, then it yields the following: If $m_n$ denotes the number of elements in $\mathcal{S}_n(\widehat{\varphi}_n) \cap T$, then for any successive knot points $t_{i-1}$ and $t_i$ in $\mathcal{S}_n(\widehat{\varphi}_n) \cap T$,

$$\sup_{i=2,\ldots,m_n} (t_i - t_{i-1}) = O_p(\rho_n^{1/5}), \tag{4.8}$$

where $\rho_n = \log(n)/n$.

PROOF OF THEOREM 4.3. From the first characterization of the estimator $\widehat{f}_n$ in Dümbgen and Rufibach (2009), for every function $\Delta$ such that $\widehat{\varphi}_n + t\Delta$ is concave for a $t > 0$ small enough, we know that

$$\int_{\mathbb{R}} \Delta(x) \, d\mathbb{F}_n(x) \leq \int_{\mathbb{R}} \Delta(x) \, d\widehat{F}_n(x). \tag{4.9}$$

This is equivalent to

$$\int_{\mathbb{R}} \Delta(x) \, d(\mathbb{F}_n(x) - F_0(x)) \leq \int_{\mathbb{R}} \Delta(x) (\widehat{f}_n(x) - f_0(x)) \, dx. \tag{4.10}$$

Using specific indicator functions for $\Delta$, one can furthermore show that

$$\widehat{F}_n(\tau) \in [\mathbb{F}_n(\tau) - 1/n, \mathbb{F}_n(\tau)] \tag{4.11}$$

for every $\tau \in \widehat{\mathcal{S}}_n(\widehat{\varphi}_n)$ [see Rufibach (2006) and Corollary 2.5 of Dümbgen and Rufibach (2009)].

Now, the idea is to choose a particular permissible perturbation function $\Delta$ that satisfies the following two conditions:

1. $\Delta$ is "local," that is, compactly supported on $[\tau_n^-, \tau_n^+]$.
2. $\Delta$ should "filter" out the unknown error $\widehat{f}_n - f_0$.

The second requirement means that $\Delta$ should be chosen so that

$$\int_{\tau_n^-}^{\tau_n^+} \Delta(x) \, dx = 0, \qquad \int_{\tau_n^-}^{\tau_n^+} \Delta(x)(x - \overline{\tau}) \, dx = 0, \tag{4.12}$$

where $\overline{\tau} := (\tau_n^- + \tau_n^+)/2$ is the mid-point of $[\tau_n^-, \tau_n^+]$. If this is guaranteed, then the right-hand side of (4.10) in the end will only depend on the distance $\tau_n^+ - \tau_n^-$ and $f_0(x_0)$.

Define $\Delta_0$ by

$$\Delta_0(x) = (x - \tau_n^-)1_{[\tau_n^-, \overline{\tau}]}(x) + (\tau_n^+ - x)1_{[\overline{\tau}, \tau_n^+]}(x).$$

Since $\widehat{\varphi}_n + t\Delta_0$ is concave for small $t > 0$, $\Delta_0$ is permissible. It is also compactly supported. However, since $\Delta_0$ is nonnegative, there is no hope that it fulfills the second of the requirements above. We therefore introduce a modified perturbation function

$$\Delta_1(x) = \Delta_0(x) - \tfrac{1}{4}(\tau_n^+ - \tau_n^-)1_{[\tau_n^-, \tau_n^+]}(x), \qquad x \in \mathbb{R}.$$

Clearly, existence of a $t > 0$, such that $\widehat{\varphi}_n + t\Delta_1$ is concave, is no longer guaranteed. However, using (4.11),

$$\int \Delta_1(x) \, d(\mathbb{F}_n - F_0)(x)$$

$$= \int \Delta_1(x) \, d(\mathbb{F}_n - \widehat{F}_n)(x) + \int \Delta_1(x) \, d(\widehat{F}_n - F_0)(x)$$

(4.13) $\qquad \leq \dfrac{\tau_n^+ - \tau_n^-}{4}\left|\int_{\tau_n^-}^{\tau_n^+} d(\mathbb{F}_n - \widehat{F}_n)(x)\right| + \int \Delta_1(x) \, d(\widehat{F}_n - F_0)(x)$

(4.14) $\qquad \leq \dfrac{\tau_n^+ - \tau_n^-}{2n} + \int \Delta_1(x)(\widehat{f}_n - f_0)(x) \, dx.$

To get the inequality in (4.13), we used (4.9) with $\Delta = \Delta_0$ and (4.11). The next step is to get bounds for the integrals in the crucial inequality (4.14). Define

$$R_{1n} := \int \Delta_1(x)(\widehat{f}_n - f_0)(x) \, dx$$

and

$$R_{2n} := \int \Delta_1(x) \, d(\mathbb{F}_n - F_0)(x).$$

Rearranging the inequality in (4.14) and using these definitions yields

$$-R_{1n} \leq \dfrac{\tau_n^+ - \tau_n^-}{2n} - R_{2n}.$$

Consistency of $\widehat{\varphi}_n$, together with $\varphi_0^{(k)}(x_0) < 0$, implies $\tau_n^+ - \tau_n^- = o_p(1)$. Thus, it follows from Lemma 4.4 that

$$M_k(-\varphi_0^{(k)}(x_0))(\tau_n^+ - \tau_n^-)^{k+2}(1 + o_p(1)) \leq o_p(1)n^{-1} + O_p(r_n^{-1}) = O_p(r_n^{-1}).$$

This yields the claimed rate, $O_p(n^{-1/(2k+1)})$, for the distance between $\tau_n^+$ and $\tau_n^-$. $\square$

LEMMA 4.4. *Suppose* (A1)–(A4) *hold. Then,*

$$R_{2n} = O_p(r_n^{-1})$$

*and*

$$R_{1n} = M_k f_0(x_0) \varphi_0^{(k)}(x_0)(\tau_n^+ - \tau_n^-)^{k+2} + o_p((\tau_n^+ - \tau_n^-)^{k+2}),$$

*where $M_k > 0$ depends only on $k$ and $\varphi_0^{(k)}(x_0) < 0$.*

PROOF. Define the function $p_n(t) = \widehat{\varphi}_n(t) - \varphi_0(t)$ for any $t \in [\tau_n^-, \tau_n^+]$. Then, using Taylor expansion of $h \mapsto \exp(h)$ up to order $k$, we can find $\theta_{t,n} \in [\tau_n^-, \tau_n^+]$, such that

$$R_{1n} = \int_{\tau_n^-}^{\tau_n^+} \Delta_1(t) f_0(t) \left( \sum_{j=1}^{k-1} \frac{p_n(t)^j}{j!} + \frac{1}{k!} \exp(\theta_{t,n}) p_n(t)^k \right) dt := \sum_{j=1}^{k} \frac{S_{nj}}{j!},$$

where

$$S_{nj} := \int_{\tau_n^-}^{\tau_n^+} \Delta_1(t) f_0(t) p_n(t)^j \, dt \qquad \text{for } 1 \leq j \leq k-1$$

and

$$S_{nk} := \int_{\tau_n^-}^{\tau_n^+} \Delta_1(t) f_0(t) \exp(\theta_{t,n}) p_n(t)^k \, dt.$$

If we expand $f_0(t)$ around the mid-point $\bar{\tau}$ of $[\tau_n^-, \tau_n^+]$, we get, for $1 \leq j \leq k-1$ and a $\eta_{n,t,j} \in [\tau_n^-, \tau_n^+]$,

$$S_{nj} = \sum_{l=0}^{k-1} \frac{f_0^{(l)}(\bar{\tau})}{l!} \int_{\tau_n^-}^{\tau_n^+} \Delta_1(t)(t-\bar{\tau})^l p_n(t)^j \, dt$$

$$+ \int_{\tau_n^-}^{\tau_n^+} \frac{f_0^{(k)}(\eta_{n,t,j})}{k!} \Delta_1(t)(t-\bar{\tau})^k p_n(t)^j \, dt$$

and, for $j = k$

$$S_{nk} = \sum_{l=0}^{k-1} \frac{f_0^{(l)}(\bar{\tau})}{l!} \int_{\tau_n^-}^{\tau_n^+} \Delta_1(t) \exp(\theta_{t,n})(t-\bar{\tau})^l p_n(t)^k \, dt$$

$$+ \int_{\tau_n^-}^{\tau_n^+} \frac{f_0^{(k)}(\eta_{n,t,k})}{k!} \Delta_1(t) \exp(\theta_{t,n})(t-\bar{\tau})^k p_n(t)^k \, dt.$$

It turns out that the dominating term in $R_{1n}$ is the first term in the Taylor expansion of $S_{n1}$. All the other terms are of smaller order since both $p_n$ and $(t-\bar{\tau})^l, l > 0$, are $o_p(1)$ uniformly in $t \in [\tau_n^-, \tau_n^+]$. We denote this dominating term by $Q_{n1}$. Since $\widehat{\varphi}_n$ is linear on $[\tau_n^-, \tau_n^+]$, we write $\widehat{\varphi}_n(t) = \widehat{\varphi}_n(\bar{\tau}) + (t - \bar{\tau})\widehat{\varphi}_n'(\bar{\tau})$. By Taylor expansion of $p_n$ around $\bar{\tau}$, we get

$$\frac{Q_{1n}}{f_0(\bar{\tau})} = \int_{\tau_n^-}^{\tau_n^+} \Delta_1(t) p_n(t)\, dt$$

$$= p_n(\bar{\tau}) \int_{\tau_n^-}^{\tau_n^+} \Delta_1(t)\, dt + p_n'(\bar{\tau}) \int_{\tau_n^-}^{\tau_n^+} \Delta_1(t)(t - \bar{\tau})\, dt$$

$$- \sum_{j=2}^{k} \frac{\varphi_0^{(j)}(\bar{\tau})}{j!} \int_{\tau_n^-}^{\tau_n^+} \Delta_1(t)(t - \bar{\tau})^j\, dt - \int_{\tau_n^-}^{\tau_n^+} \varepsilon_n(t) \Delta_1(t)(t - \bar{\tau})^k\, dt,$$

where the first two terms are zero, since (4.12) holds when $\Delta = \Delta_1$ and $\|\varepsilon_n\|_\infty \to_p 0$ as $\tau_n^+ - \tau_n^- \to_p 0$. Using the fact that

(4.15)
$$\int_{\tau_n^-}^{\tau_n^+} \Delta_1(t)(t-\bar{\tau})^j\, dt$$
$$= \begin{cases} 0, & \text{for } j = 0 \text{ and } j \text{ odd,} \\ (\tau_n^+ - \tau_n^-)^{j+2}\left(\dfrac{-j}{2^{(j+2)}(j+1)(j+2)}\right), & \text{for } j \text{ even,} \end{cases}$$

we conclude that

$$Q_{1n} = \frac{k}{2^{(k+2)} k!(k+1)(k+2)} f_0(\bar{\tau}) \varphi_0^{(k)}(\bar{\tau})((\tau_n^+ - \tau_n^-)^{k+2} + o_p(1))$$

and the claimed form of $R_{1n}$ in the lemma follows.

For $R_{2n}$, we proceed along the lines of the proof of Lemma 4.1 in Groeneboom, Jongbloed and Wellner (2001b). This means we have to line up with the assumption of Theorem 2.14.1 in van der Vaart and Wellner (1996). Therefore, define a generalized version of $R_{2n}$:

$$R_{2n}^{x,y} = \int_x^y \Delta_1(z)\, d(\mathbb{F}_n - F_0)(z)$$

for $-\infty < x \leq y$. With this function, we have, for some $R > 0$,

$$\sup_{y\,:\,0 \leq y-x \leq R} |R_{2n}^{x,y}|$$

$$= 2 \sup_{y\,:\,0 \leq y-x \leq R} \left| \int_x^{(x+y)/2} (z - x - \tfrac{1}{4}(y-x))\, d(\mathbb{F}_n - F_0)(z) \right|$$

$$= 2 \sup_{y\,:\,0 \leq y-x \leq R} \left| \int h_{x,y}(z)\, d(\mathbb{F}_n - F_0)(z) \right|,$$

where

$$h_{x,y}(z) = (z - x - \tfrac{1}{4}(y-x)) 1_{[x,(x+y)/2]}(z) = h(z) 1_{[x,(x+y)/2]}(z).$$

Then, the collection of functions

$$\mathcal{F}_{x,R} = \{h 1_{[x,(x+y)/2]} : x \leq y \leq x+R\}$$

is a Vapnik–Chervonenkis subgraph class with envelope function

$$F_{x,R}(z) = ((z-x) + R/4) 1_{[x,x+R]}(z).$$

Finally, Theorem 2.6.7 in van der Vaart and Wellner (1996) yields the entropy condition (4.5).

A log-concave density is always unimodal and the value at the mode is finite, and hence, $K := \|f_0\|_\infty$ is finite. Therefore,

$$EF_{x,R}^2(X_1)$$
$$= \int_x^{x+R} (z-x)^2 f_0(z)\, dz + \frac{R}{2} \int_x^{x+R} (z-x) f_0(z)\, dz + \frac{R^2}{16} \int_x^{x+R} f_0(z)\, dz$$
$$\leq \left( \frac{K}{3}(z-x)^3 + \frac{RK}{4}(z-x)^2 + \frac{R^2 K}{16} z \right) \Big|_{z=x}^{x+R}$$
$$= \frac{31}{48} K R^3.$$

It follows from Lemma 4.1, with $d=2$ and $s=k$, that $R_{2n} = O_p(r_n^{-1})$. □

4.2. *Proofs for Section 2.*

LEMMA 4.5. *For any $M > 0$, we have*

$$(4.16) \qquad \sup_{|t| \leq M} |\widehat{\varphi}_n'(x_0 + s_n t) - \varphi_0'(x_0)| = O_p(s_n^{k-1}),$$

$$(4.17) \qquad \sup_{|t| \leq M} |\widehat{\varphi}_n(x_0 + s_n t) - \varphi_0(x_0) - s_n t \varphi_0'(x_0)| = O_p(s_n^k).$$

*Furthermore, if we define, for any $u \in \mathbb{R}$,*

$$\hat{e}_n(u) = \widehat{f}_n(u) - \sum_{j=0}^{k-1} \frac{f_0^{(j)}(x_0)}{j!} (u - x_0)^j - f_0(x_0) \frac{[\varphi_0'(x_0)]^k}{k!} (u - x_0)^k,$$

*then*

$$(4.18) \qquad \sup_{|t| \leq M} |\hat{e}_n(x_0 + s_n t) - f_0(x_0)(\widehat{\varphi}_n(x_0 + s_n t) - \varphi_0(x_0) - s_n t \varphi_0'(x_0))|$$
$$= o_p(s_n^k).$$

PROOF. The proof of (4.16) and (4.17) is identical to that of Lemma 4.4 in Groeneboom, Jongbloed and Wellner (2001b) since the characterization of $\widehat{f}_n$ given in (1.1) is (up to the direction of the inequality) equivalent to that of the least-squares estimator of a convex density.

Now, we prove (4.18). Using Taylor expansion of $h \mapsto \exp(h)$ up to order $k$ around zero, we can write

$$
\begin{aligned}
(4.19) \quad \widehat{f}_n(u) - f_0(x_0) &= f_0(x_0)[\exp(\widehat{\varphi}_n(u) - \varphi_0(x_0)) - 1] \\
&= f_0(x_0) \sum_{j=1}^{k} \frac{1}{j!}(\widehat{\varphi}_n(u) - \varphi_0(x_0))^j + f_0(x_0)\widehat{\Psi}_{k,n,1}(u),
\end{aligned}
$$

where

$$\widehat{\Psi}_{k,n,1}(u) = \sum_{j=k+1}^{\infty} \frac{1}{j!}(\widehat{\varphi}_n(u) - \varphi_0(x_0))^j.$$

But, for any $j \geq 1$,

$$
\begin{aligned}
(\widehat{\varphi}_n(u) - \varphi_0(x_0))^j &= [\widehat{\varphi}_n(u) - \varphi_0(x_0) - (u-x_0)\varphi_0'(x_0) + (u-x_0)\varphi_0'(x_0)]^j \\
(4.20) \quad &= \sum_{r=1}^{j} \binom{j}{r} [\widehat{\varphi}_n(u) - \varphi_0(x_0) - (u-x_0)\varphi'(x_0)]^r \\
&\quad \times [\varphi_0'(x_0)]^{j-r}(u-x_0)^{j-r} \\
&\quad + [\varphi_0'(x_0)]^j (u-x_0)^j.
\end{aligned}
$$

Hence, using (4.17) and (A3), we get on the set $\{u : |u-x_0| \leq Mn^{-1/(2k+1)}\}$

$$(\widehat{\varphi}_n(u) - \varphi_0(x_0))^j = o_p(n^{-k/(2k+1)})$$

for all $j \geq k+1$.

In particular, this implies that

$$(4.21) \quad \widehat{\Psi}_{k,n,1}(u) = o_p(n^{-k/(2k+1)}),$$

uniformly in $u \in [x_0 - tn^{-1/(2k+1)}, x_0 + tn^{-1/(2k+1)}]$, where $|t| \leq M$, and

$$
\begin{aligned}
&\widehat{f}_n(u) - f_0(x_0) - f_0(x_0)(\widehat{\varphi}_n(u) - \varphi_0(x_0) - (u-x_0)\varphi_0'(x_0)) \\
&\quad - f_0(x_0)\sum_{j=1}^{k} \frac{\varphi_0^{(j)}(x_0)}{j!}(u-x_0)^j = o_p(n^{-k/(2k+1)}).
\end{aligned}
$$

Using Lemma 4.2, the latter can be rewritten as

$$
\begin{aligned}
&\widehat{f}_n(u) - f_0(x_0) - f_0(x_0)(\widehat{\varphi}_n(u) - \varphi_0(x_0) - (u-x_0)\varphi_0'(x_0)) \\
&\quad - \sum_{j=1}^{k-1} \frac{f_0^{(j)}(x_0)}{j!}(u-x_0)^j - f_0(x_0)\frac{\varphi_0^{(k)}(x_0)}{k!}(u-x_0)^k = o_p(n^{-k/(2k+1)})
\end{aligned}
$$

or, equivalently,

$$|\hat{e}_n(x_0 + tn^{-1/(2k+1)}) - f_0(x_0)(\widehat{\varphi}_n(x_0 + tn^{-1/(2k+1)})$$
$$- \varphi_0(x_0) - n^{-1/(2k+1)}t\varphi_0'(x_0))| = o_p(n^{-k/(2k+1)})$$

uniformly in $|t| \leq M$.  $\square$

THEOREM 4.6. *Let $K > 0$.*

(i) *If $\{Y_k(t), t \in \mathbb{R}\}$ is the canonical process defined in (2.1), then the localized process $\gamma_1 \mathbb{Y}_n^{\mathrm{locmod}}(\gamma_2 \cdot)$ converges weakly in $C[-K, K]$ to $Y_k$, where*

$$(4.22) \qquad \gamma_1 = \left( \frac{f_0(x_0)^{k-1}|\varphi_0^{(k)}(x_0)|^3}{[(k+2)!]^3} \right)^{1/(2k+1)},$$

$$(4.23) \qquad \gamma_2 = \left( \frac{f_0(x_0)|\varphi_0^{(k)}(x_0)|^2}{[(k+2)!]^2} \right)^{1/(2k+1)}.$$

*Equivalently, $\mathbb{Y}_n^{\mathrm{locmod}}$ converges weakly in $C[-K, K]$ to the "driving process" $Y_{a,k,\sigma}$, where*

$$(4.24) \qquad Y_{k,a,\sigma}(t) := a \int_0^t W(s)\, ds - \sigma t^{k+2}$$

*and where $a = 1/\sqrt{f_0(x_0)}$, $\sigma = |\varphi_0^{(k)}(x_0)|/(k+2)!$.*

(ii) *The localized processes satisfy $\mathbb{Y}_n^{\mathrm{locmod}}(t) - \widehat{H}_n^{\mathrm{locmod}}(t) \geq 0$, for all $t \in \mathbb{R}$, with equality for all $t$ such that $x_n(t) = x_0 + tn^{-1/(2k+1)} \in \widehat{\mathcal{S}}_n(\widehat{\varphi}_n)$.*

(iii) *Both $\hat{A}_n$ and $\hat{B}_n$ defined above in (4.1) and (4.2) are tight.*

(iv) *The vector of processes*

$$(\widehat{H}_n^{\mathrm{locmod}}, (\widehat{H}_n^{\mathrm{locmod}})^{(1)}, (\widehat{H}_n^{\mathrm{locmod}})^{(2)}, \mathbb{Y}_n^{\mathrm{locmod}}, (\widehat{H}_n^{\mathrm{locmod}})^{(3)}, (\mathbb{Y}_n^{\mathrm{locmod}})^{(1)})$$

*converges weakly in $(C[-K,K])^4 \times (D[-K,K])^2$, endowed with the product topology induced by the uniform topology on the spaces $C[-K,K]$ and the Skorohod topology on the spaces $D[-K,K]$ to the process*

$$(H_{k,a,\sigma}, H_{k,a,\sigma}^{(1)}, H_{k,a,\sigma}^{(2)}, Y_{k,a,\sigma}, H_{k,a,\sigma}^{(3)}, Y_{k,a,\sigma}^{(1)}),$$

*where $H_{k,a,\sigma}$ is the unique process on $\mathbb{R}$ satisfying*

$$(4.25) \quad \begin{cases} H_{k,a,\sigma}(t) \leq Y_{k,a,\sigma}(t), & \text{for all } t \in \mathbb{R}, \\ \int (H_{k,a,\sigma}(t) - Y_{k,a,\sigma}(t))\, dH_{k,a,\sigma}^{(3)}(t) = 0, \\ H_{k,a,\sigma}^{(2)}, & \text{is concave}. \end{cases}$$

PROOF. (i) The first step will be to modify the local processes, that is, going from the "density" to the "log-density" level, in order to be able to exploit concavity of $\varphi_0$ and $\widehat{\varphi}_n$ and connect the local process to the limiting distribution obtained by Groeneboom, Jongbloed and Wellner (2001b) for estimating a convex density.

First, by Lemma 4.2, (4.19) and (A3), we can write

$$f_0(x_0)^{-1}\left(\widehat{f}_n(u) - \sum_{j=0}^{k-1} \frac{f_0^{(j)}(x_0)}{j!}(u-x_0)^j\right)$$

$$= f_0(x_0)^{-1}\left(\widehat{f}_n(u) - f_0(x_0) - f_0(x_0)\sum_{j=1}^{k-1} \frac{[\varphi_0'(x_0)]^j}{j!}(u-x_0)^j\right)$$

$$= \widehat{\Psi}_{k,n,1}(u) + \sum_{j=1}^{k} \frac{1}{j!}[\widehat{\varphi}_n(u) - \varphi_0(x_0)]^j - \sum_{j=1}^{k-1} \frac{[\varphi_0'(x_0)]^j}{j!}(u-x_0)^j$$

$$= \widehat{\Psi}_{k,n,1}(u) + (\widehat{\varphi}_n(u) - \varphi_0(x_0) - \varphi_0'(x_0)(u-x_0))$$

$$+ \sum_{j=2}^{k} \frac{1}{j!}[\widehat{\varphi}_n(u) - \varphi_0(x_0)]^j - \sum_{j=2}^{k-1} \frac{[\varphi_0'(x_0)]^j}{j!}(u-x_0)^j$$

$$=: (\widehat{\varphi}_n(u) - \varphi_0(x_0) - \varphi_0'(x_0)(u-x_0)) + \widehat{\Psi}_{k,n,2}(u),$$

introducing the new remainder term

(4.26)
$$\widehat{\Psi}_{k,n,2}(u) = \widehat{\Psi}_{k,n,1}(u) + \sum_{j=2}^{k} \frac{1}{j!}[\widehat{\varphi}_n(u) - \varphi_0(x_0)]^j$$
$$- \sum_{j=2}^{k-1} \frac{[\varphi_0'(x_0)]^j}{j!}(u-x_0)^j.$$

Using (4.20) and (4.21) yields

$$\int_{\mathcal{I}} \int_{x_0}^{v} \widehat{\Psi}_{k,n,2}(u)\,du\,dv$$

$$= t^2 n^{-2/(2k+1)} \sup_{u\in[x_0,v], v\in\mathcal{I}} |\widehat{\Psi}_{k,n,1}(u)|$$

$$+ \sum_{j=2}^{k} \frac{1}{j!} \int_{\mathcal{I}} \int_{x_0}^{v} [\widehat{\varphi}_n(u) - \varphi_0(x_0)]^j\,du\,dv$$

$$- \sum_{j=2}^{k-1} \frac{1}{j!} \int_{\mathcal{I}} \int_{x_0}^{v} [\varphi_0'(x_0)]^j(u-x_0)^j\,du\,dv$$

$$= o_p(r_n^{-1})$$
$$+ \sum_{j=2}^{k} \frac{1}{j!} \sum_{l=1}^{j} \binom{j}{l} \int_{\mathcal{I}} \int_{x_0}^{v} [\widehat{\varphi}_n(u) - \varphi_0(x_0)$$
$$- (u-x_0)\varphi_0'(x_0)]^l$$
$$\times (u-x_0)^{j-l} [\varphi_0'(x_0)]^{j-l} \, du \, dv$$
$$+ \sum_{j=2}^{k} \frac{1}{j!} \int_{\mathcal{I}} \int_{x_0}^{v} [\varphi_0'(x_0)]^j (u-x_0)^j \, du \, dv$$
$$- \sum_{j=2}^{k-1} \frac{1}{j!} \int_{\mathcal{I}} \int_{x_0}^{v} [\varphi_0'(x_0)]^j (u-x_0)^j \, du \, dv$$
$$= o_p(r_n^{-1})$$
$$+ \sum_{j=2}^{k} \frac{1}{j!} \sum_{l=1}^{j} \binom{j}{l} \int_{\mathcal{I}} \int_{x_0}^{v} [\widehat{\varphi}_n(u) - \varphi_0(x_0)$$
$$- (u-x_0)\varphi_0'(x_0)]^l$$
$$\times (u-x_0)^{j-l} [\varphi_0'(x_0)]^{j-l} \, du \, dv$$
$$+ \frac{1}{k!} \int_{\mathcal{I}} \int_{x_0}^{v} (u-x_0)^k [\varphi_0'(x_0)]^k \, du \, dv.$$

But by Lemma 4.5, one can easily show that, for $j = 2, \ldots, k$ and $l = 1, \ldots, j$,

$$r_n \int_{\mathcal{I}} \int_{x_0}^{v} [\widehat{\varphi}_n(u) - \varphi_0(x_0) - (u-x_0)\varphi_0'(x_0)]^l (u-x_0)^{j-l} [\varphi_0'(x_0)]^{j-l} \, du \, dv$$
$$= O_p(n^{-[k(l-1)+(j-l)]/(2k+1)}) = o_p(1),$$

uniformly in $|t| \leq M$. Similarly,

$$r_n \int_{\mathcal{I}} \int_{x_0}^{v} (u-x_0)^k [\varphi_0'(x_0)]^k \, du \, dv = \frac{[\varphi_0'(x_0)]^k}{(k+1)(k+2)} t^{k+2}.$$

Hence, it follows that

$$r_n \int_{\mathcal{I}} \int_{x_0}^{v} \widehat{\Psi}_{k,n,2}(u) \, du \, dv = \frac{[\varphi_0'(x_0)]^k}{(k+2)!} t^{k+2} + o_p(1)$$

as $n \to \infty$, uniformly in $|t| \leq M$.

We turn now to the modified local processes, $\mathbb{Y}_n^{\text{locmod}}$ and $\widehat{H}_n^{\text{locmod}}$, defined in (4.3) and (4.4). It is not difficult to show that

(4.27) $$\mathbb{Y}_n^{\text{locmod}}(t) = \frac{\mathbb{Y}_n^{\text{loc}}(t)}{f_0(x_0)} - r_n \int_{\mathcal{I}} \int_{x_0}^{v} \widehat{\Psi}_{k,n,2}(u) \, du \, dv$$

and

$$(4.28) \quad \widehat{H}_n^{\text{locmod}}(t) = \frac{\widehat{H}_n^{\text{loc}}(t)}{f_0(x_0)} - r_n \int_{\mathcal{I}} \int_{x_0}^{v} \widehat{\Psi}_{k,n,2}(u) \, du \, dv.$$

Note that the process $\widehat{H}_n^{\text{locmod}}$ is in fact similar to $\widehat{H}_n^{\text{loc}}$, except that it is defined in terms of the log-density $\varphi_0$ instead of the density $f_0$. This can be more easily seen from its original expression given in (4.4). The second expression of $\widehat{H}_n^{\text{locmod}}$ given above is only useful for showing that it stays below $\mathbb{Y}_n^{\text{locmod}}$, while touching it at points $t$, such that $x_n(t) = x_0 + tn^{-1/(2k+1)} \in \widehat{\mathcal{S}}_n(\widehat{\varphi}_n)$. The biggest advantage of considering this modified version is to be able to use concavity of $\varphi_0$ the same way [Groeneboom, Jongbloed and Wellner (2001b)] used convexity of the true estimated density $g_0$. Their process $\widetilde{H}_n^{\text{loc}}$ resembles $\widehat{H}_n^{\text{locmod}}$ to a large extent (see page 1688), and by combining arguments similar to theirs with Lemma 4.2 and the results obtained above, it follows that

$$\mathbb{Y}_n^{\text{locmod}}(t)$$

$$\Rightarrow [f_0(x_0)]^{-1/2} \int_0^t W(s) \, ds + \frac{f_0^{(k)}(x_0)}{(k+2)! f_0(x_0)} t^{k+2} - \frac{[\varphi_0'(x_0)]^k}{(k+2)!} t^{k+2}$$

$$= [f_0(x_0)]^{-1/2} \int_0^t W(s) \, ds + \frac{\varphi_0^{(k)}(x_0)}{(k+2)!} t^{k+2}$$

$$= Y_{k,a,\sigma}(t) \quad \text{in } C[-K,K],$$

where $a := [f_0(x_0)]^{-1/2}$, $\sigma := |\varphi_0^{(k)}(x_0)|/(k+2)!$, as in (4.24).

Now, let $\gamma_1$ and $\gamma_2$ be chosen, so that

$$\gamma_1 Y_{k,a,\sigma}(\gamma_2 t) \stackrel{d}{=} Y_k(t)$$

as processes where $Y_k$ is the integrated Gaussian process defined in (2.1). Using the scaling property of Brownian motion [i.e., $\alpha^{-1/2} W(\alpha t) \stackrel{d}{=} W(t)$, for any $\alpha > 0$], we get

$$\gamma_1 \gamma_2^{3/2} = a^{-1} \quad \text{and} \quad \gamma_1 \gamma_2^{k+2} = \sigma^{-1}.$$

This yields $\gamma_1$ and $\gamma_2$ as given in (4.22) and (4.23), and hence,

$$\begin{pmatrix} n^{k/(2k+1)}(\widehat{\varphi}_n(x_0) - \varphi_0(x_0)) \\ n^{(k-1)/(2k+1)}(\widehat{\varphi}_n'(x_0) - \varphi_0'(x_0)) \end{pmatrix} \stackrel{d}{\to} f_0(x_0)^{-1} \begin{pmatrix} c_k(x_0, \varphi_0) H_k^{(2)}(0) \\ d_k(x_0, \varphi_0) H_k^{(3)}(0) \end{pmatrix}.$$

We get the explicit expression of the asymptotic constants $c_k(x_0, \varphi_0)$ and $d_k(x_0, \varphi_0)$ using the following relations:

$$(4.29) \quad f_0(x_0)^{-1} c_k(x_0, \varphi_0) = (\gamma_1 \gamma_2^2)^{-1} \quad \text{and}$$

$$(4.30) \quad f_0(x_0)^{-1} d_k(x_0, \varphi_0) = (\gamma_1 \gamma_2^3)^{-1}.$$

This is completely analogous to the derivations on page 1689 in Groeneboom, Jongbloed and Wellner (2001b), precisely

$$(4.31) \quad \begin{aligned} (\gamma_1 \widehat{H}_n^{\text{locmod}}(\gamma_2 t))^{(2)}(0) &= \gamma_1 \gamma_2^2 (\widehat{H}_n^{\text{locmod}})^{(2)}(0) \\ &= n^{k/(2k+1)} f_0(x_0) c_k(x_0, \varphi_0)^{-1} (\widehat{\varphi}_n(x_0) - \varphi_0(x_0)) \end{aligned}$$

and

$$(4.32) \quad \begin{aligned} (\gamma_1 \widehat{H}_n^{\text{locmod}}(\gamma_2 t))^{(3)}(0) &= \gamma_1 \gamma_2^3 (\widehat{H}_n^{\text{locmod}})^{(3)}(0) \\ &= n^{(k-1)/(2k+1)} f_0(x_0) d_k(x_0, \varphi_0)^{-1} (\widehat{\varphi}_n'(x_0) - \varphi_0'(x_0)). \end{aligned}$$

From (4.29) and (4.30), we get $c_k(x_0, \varphi_0)$ and $d_k(x_0, \varphi_0)$ as given in (2.2) and (2.3), and $C_k(x_0, \varphi_0)$ and $D_k(x_0, \varphi_0)$ as in (2.4) and (2.5).

(ii) Note that we can write

$$\mathbb{Y}_n^{\text{loc}}(t) - \widehat{H}_n^{\text{loc}}(t) = r_n(\mathbb{H}_n(x_n(t)) - \widehat{H}_n(x_n(t))) \geq 0$$

by making use of (1.1) and the specific choice of $\hat{A}_n$ and $\hat{B}_n$. But, since we connect $\widehat{H}_n^{\text{locmod}}$ and $\mathbb{Y}_n^{\text{locmod}}$ to the "invelope," the latter property needs primarily to hold for the modified processes. This can easily be established by considering (4.27) and (4.28), and hence it follows that

$$\mathbb{Y}_n^{\text{locmod}}(t) - \widehat{H}_n^{\text{locmod}}(t) \geq 0$$

for all $t \in \mathbb{R}$, with equality if $x_n(t) = x_0 + t n^{-1/(2k+1)} \in \widehat{\mathcal{S}}_n(\widehat{\varphi}_n)$.

(iii) To show that $\hat{A}_n$ and $\hat{B}_n$ are tight. By Theorem 4.3, we know that there exists $M > 0$ and $\tau \in \widehat{\mathcal{S}}(\widehat{\varphi}_n)$ such that $0 \leq x_0 - \tau \leq M n^{-1/(2k+1)}$ with large probability. Now, using (4.11), we can write

$$\begin{aligned} |\hat{A}_n| &\leq r_n s_n |(\widehat{F}_n(x_0) - \widehat{F}_n(\tau)) - (\mathbb{F}_n(x_0) - \mathbb{F}_n(\tau))| + r_n/n \\ &\leq r_n s_n \left| \int_\tau^{x_0} \left( \widehat{f}_n(u) - \sum_{j=0}^{k-1} \frac{f_0^{(j)}(x_0)}{j!} (u - x_0)^j \right) du \right| \\ &\quad + r_n s_n \left| \int_\tau^{x_0} \left( \sum_{j=0}^{k-1} \frac{f_0^{(j)}(x_0)}{j!} (u - x_0)^j - f_0(u) \right) du \right| \\ &\quad + r_n s_n \left| \int_\tau^{x_0} d(\mathbb{F}_n - F_0) \right| + n^{-k/(2k+1)} \\ &:= \hat{A}_{n1} + \hat{A}_{n2} + \hat{A}_{n3} + n^{-k/(2k+1)}. \end{aligned}$$

Now,

$$|\hat{A}_{n1}| \leq r_n s_n \left| \int_\tau^{x_0} \hat{e}_n(u) \, du - f_0(x_0)(\widehat{\varphi}_n(u) - \varphi_0(x_0) - (u - x_0) \varphi_0'(x_0)) \, du \right|$$

$$+ r_n s_n f_0(x_0) \left| \int_\tau^{x_0} \left( \frac{[\varphi_0'(x_0)]^k}{k!} (u - x_0)^k \right) du \right|$$

$$+ r_n s_n f_0(x_0) \left| \int_\tau^{x_0} (\widehat{\varphi}_n(u) - \varphi_0(x_0) - (u - x_0)\varphi_0'(x_0))\, du \right|$$

$$\leq o_p(1) + O_p(r_n s_n (\tau - x_0)^{k+1}) + O_p(r_n s_n (\tau - x_0) n^{-k/(2k+1)})$$

$$= O_p(1),$$

where we used (4.18) and (4.17) to bound the first and last terms. To bound $\hat{A}_{n2}$, we use Taylor approximation of $f_0(u)$ around $x_0$ to get

$$\hat{A}_{n2} \leq r_n \left| \int_\tau^{x_0} \frac{f_0^{(k)}(x_0)}{k!} (u - x_0)^k \, du \right| + r_n \left| \int_\tau^{x_0} (u - x_0)^k \varepsilon_n(u) \, du \right|$$

$$= O_p(1),$$

where $\varepsilon_n$ is a function such that $\|\varepsilon_n\| \to_p 0$ as $x_0 - \tau \to_p 0$. To bound $\hat{A}_{n3}$, similar derivations as the ones used for bounding $R_{2n}$ (see the proof of Lemma 4.4) can be employed where the perturbation function $\Delta_1$ needs to be replaced by $\Delta_2(x) = 1_{[\tau, x_0]}(x)$.

At "one higher integration level," similar computations can be used to show tightness of $\hat{B}_n$.

(iv) The proof of this last part of the theorem is basically identical to that of Theorem 6.2 for the LSE in Groeneboom, Jongbloed and Wellner (2001b) and arguments similar to those of Groeneboom, Jongbloed and Wellner (2001a) or, alternatively, tightness plus uniqueness arguments along the lines of Groeneboom, Maathuis, and Wellner (2008). □

PROOF OF THEOREM 2.1. The claimed joint convergence involving $\widehat{\varphi}_n$ and $\widehat{\varphi}_n'$ follows from part (iv) of Theorem 4.6 and the relations (4.31) and (4.32). The joint limiting distribution of $\widehat{f}_n(x_0) - f_0(x_0)$ and $\widehat{f}_n'(x_0) - f_0'(x_0)$ follows immediately by applying the delta-method. □

4.3. *Proofs for Section 3.*

PROOF OF THEOREM 3.6. We first use the simple fact that $\widehat{M}_n$ is the only point $x \in \mathbb{R}$ which satisfies

(4.33) $$\widehat{\varphi}_n'(t) \begin{cases} > 0, & \text{if } t < x, \\ \leq 0, & \text{if } t \geq x. \end{cases}$$

This follows immediately from concavity of $\widehat{\varphi}_n$ and the definition of $\widehat{M}_n$. Note that $\widehat{\varphi}_n$ may have a flat region or "modal interval"; in this case, there exists an entire interval of points where the maximum is attained, and $\widehat{M}_n$ is the left endpoint of this interval.

A tightness property of the process $H_2^{(3)}$, which follows from Lemma 2.7 of Groeneboom, Jongbloed and Wellner (2001b), is also needed to establish the limiting distribution of $\widehat{M}_n$: for any $\varepsilon > 0$ and $t \in \mathbb{R}$, there exists $C = C(\varepsilon)$ such that

$$P(|H_2^{(3)}(t) + 24t| > C) \leq \varepsilon.$$

In other words, one can view $H_2^{(3)}(t)$ as an "estimator" of the odd function $-24t$. Since $C$ is independent of $t$, it follows that, for a fixed $\varepsilon$, $H_2^{(3)}(t) < 0$ [resp. $H_2^{(3)}(t) > 0$] for $t > 0$ (resp. $-t < 0$) big enough, with probability greater than $1 - \varepsilon$.

The sign of $H_2^{(3)}$ and uniqueness of $\widehat{M}_n$ turn out to be crucial in determining the limiting distribution of the latter. From Theorem 4.6 and the two derivative relations, (4.31) and (4.32), it follows that

$$\begin{aligned}
(4.34) \quad & \begin{pmatrix} n^{k/(2k+1)}(\hat{\varphi}_n(x_0 + tn^{-1/(2k+1)}) - \varphi_0(x_0) - tn^{-1/(2k+1)}\varphi_0'(x_0)) \\ n^{(k-1)/(2k+1)}(\hat{\varphi}_n'(x_0 + tn^{-1/(2k+1)}) - \varphi_0'(x_0)) \end{pmatrix} \\
& \Rightarrow \begin{pmatrix} H_{k,a,\sigma}^{(2)}(t) \\ H_{k,a,\sigma}^{(3)}(t) \end{pmatrix} \quad \text{in } C[-K,K] \times D[-K,K]
\end{aligned}$$

for each $K > 0$, with the product topology induced by the uniform topology on $C[-K,K]$ and the Skorohod topology on $D[-K,K]$. Here, $H_{k,a,\sigma}$ is the unique process on $\mathbb{R}$ satisfying (4.25). A similar result holds for the MLE of the log-concave density $f_0$. When $x_0$ is replaced by the population mode $m_0 = M(f_0)$ and $k = 2$ the second weak convergence implies that

$$n^{1/5}(\hat{\varphi}_n'(m_0 - Tn^{-1/5}) - \varphi_0'(m_0)) \xrightarrow{d} H_{2,a,\sigma}^{(3)}(-T)$$

and

$$n^{1/5}(\hat{\varphi}_n'(m_0 + Tn^{-1/5}) - \varphi_0'(m_0)) \xrightarrow{d} H_{2,a,\sigma}^{(3)}(T).$$

For $T > 0$ large enough, this in turn implies that, for $\varepsilon > 0$, we can find $N \in \mathbb{N} \setminus \{0\}$ such that, for all $n > N$, we have

$$P(\hat{\varphi}_n'(m_0 - Tn^{-1/5}) > 0 \text{ and } \hat{\varphi}_n'(m_0 + Tn^{-1/5}) < 0) > 1 - \varepsilon.$$

Using the property of $\widehat{M}_n$ in (4.33), it follows that

$$P(\widehat{M}_n \in [m_0 - Tn^{-1/5}, m_0 + Tn^{-1/5}]) > 1 - \varepsilon$$

for all $n > N$.

We first conclude that $\widehat{M}_n - m_0 = O_p(n^{-1/5})$. Then, we note that

$$n^{1/5}(\widehat{M}_n - m_0) = M(\mathbb{Z}_n),$$

where

$$\mathbb{Z}_n(t) = n^{2/5}(\hat{\varphi}_n(m_0 + tn^{-1/5}) - \varphi_0(m_0))$$
$$\Rightarrow \mathbb{Z}(t) := H_{2,a,\sigma}^{(2)}(t) \qquad \text{in } C([-K,K])$$

for each $K > 0$, by (4.34) with $k = 2$. Thus, by the argmax continuous mapping theorem [see, e.g., van der Vaart and Wellner (1996), page 286] it follows that

$$M(\mathbb{Z}_n) \xrightarrow{d} M(\mathbb{Z}) = M(H_{2,a,\sigma}^{(2)}),$$

where $\mathbb{Z} = H_{2,a,\sigma}^{(2)}$, $a = 1/\sqrt{f_0(m_0)}$, and $\sigma = |\varphi_0^{(2)}(m_0)|/4!$.

Note that $H_{2,a,\sigma}$ is related to the "driving process" $Y_{2,a,\sigma}$ with $a = 1/\sqrt{f_0(m_0)}$, $\sigma = |\varphi_0^{(2)}(m_0)|/4!$ as in (4.24) with $k = 2$. Now, $\gamma_1 Y_{2,a,\sigma}(\gamma_2 t) \stackrel{d}{=} Y_2(t)$ as processes where $Y_2 := Y_{2,1,1}$. Thus, it also holds that

$$\gamma_1 H_{2,a,\sigma}(\gamma_2 t) \stackrel{d}{=} H_2(t) \quad \text{and} \quad \gamma_1 \gamma_2^2 H_{2,a,\sigma}^{(2)}(\gamma_2 t) \stackrel{d}{=} H_2^{(2)}(t),$$

or, equivalently, $H_{2,a,\sigma}^{(2)}(v) \stackrel{d}{=} H_2^{(2)}(v/\gamma_2)/(\gamma_1 \gamma_2^2)$. Since $M(dg(c\cdot)) = c^{-1}M(g)$ for $c, d > 0$, it follows that

$$M(H_{2,a,\sigma}^{(2)}) \stackrel{d}{=} M\left(\frac{1}{\gamma_1 \gamma_2^2} H_2^{(2)}(\cdot/\gamma_2)\right) \stackrel{d}{=} \gamma_2 M(H_2^{(2)}),$$

where

$$\gamma_2 = \left(f_0(m_0)\frac{|\varphi_0^{(2)}(m_0)|^2}{(4!)^2}\right)^{-1/5} = \left(\frac{(4!)^2 f_0(m_0)}{f_0''(m_0)^2}\right)^{1/5}$$

by direct computation using $f_0'(m_0) = 0 = \varphi_0'(m_0)$ and Lemma 4.2. □

**Acknowledgments.** This research was initiated while Kaspar Rufibach was visiting the Institute for Mathematical Stochastics at the University of Göttingen, Germany and while Fadoua Balabdaoui was visiting the Institute of Mathematical Statistics and Actuarial Science at the University of Bern, Switzerland. We would like to thank both institutions for their hospitality. We also thank two referees and an Associate Editor for several suggestions leading to improvements of the exposition.

F. Balabdaoui  
Centre de Recherche en Mathematiques de la Decision  
Universite Paris-Dauphine  
Paris  
France  
E-mail: fadoua@ceremade.dauphine.fr  

K. Rufibach  
Institute of Social  
and Preventive Medicine  
Biostatistics Unit  
University of Zurich  
CH-8001 Zurich  
Switzerland  
E-mail: kaspar.rufibach@ifspm.uzh.ch  

J. A. Wellner  
Department of Statistics  
University of Washington  
Box 354322  
Seattle, Washington 98195-4322  
USA  
E-mail: jaw@stat.washington.edu